\documentclass[12pt]{article}
\usepackage{amssymb, epsfig, amssymb, latexsym, color,graphicx}
\begin{document}

\newcommand{\s}{\sigma}
\renewcommand{\k}{\kappa}
\newcommand{\p}{\partial}
\newcommand{\D}{\Delta}
\newcommand{\om}{\omega}
\newcommand{\Om}{\Omega}
\renewcommand{\phi}{\varphi}
\newcommand{\e}{\epsilon}
\renewcommand{\a}{\alpha}
\renewcommand{\b}{\beta}
\newcommand{\N}{{\mathbb N}}
\newcommand{\R}{{\mathbb R}}
   \newcommand{\eps}{\varepsilon}
   \newcommand{\EX}{{\Bbb{E}}}
   \newcommand{\PX}{{\Bbb{P}}}

\newcommand{\cF}{{\cal F}}
\newcommand{\cG}{{\cal G}}
\newcommand{\cD}{{\cal D}}
\newcommand{\cO}{{\cal O}}

\newtheorem{theorem}{Theorem}
\newtheorem{lemma}{Lemma}
\newtheorem{remark}{Remark}
\newtheorem{example}{Example}
\newtheorem{definition}{Definition}

\title{Invariant manifold reduction\\ for  stochastic dynamical systems }

\author{Aijun Du and Jinqiao Duan \\
Department of Applied Mathematics \\
Illinois Institute of Technology \\
Chicago, IL 60616, USA.  $\;\;$  E-mail:  \emph{duan@iit.edu}   }

 \date{February 28, 2006   }

\maketitle

\begin{abstract}

Invariant manifolds facilitate the understanding of nonlinear
stochastic dynamics. When an invariant manifold is represented
approximately  by   a graph for example,  the whole stochastic
dynamical system may be reduced or restricted to this manifold.
This reduced system may provide valuable dynamical information for
the original system. The authors have derived an invariant
manifold reduction or restriction principle for systems of
Stratonovich or Ito stochastic differential equations.

Two concepts of invariance are considered for invariant manifolds.

The first invariance   concept   is   in the framework of cocycles
--- an invariant manifold being a random set. The dynamical reduction is
achieved by investigating random center manifolds.

The second invariance concept is in the sense of almost sure ---
an invariant manifold being a deterministic set which is not
necessarily   attracting. The   restriction of the original
stochastic system on this deterministic local invariant manifold
is still a stochastic system but with reduced dimension.





\end{abstract}



\section{Introduction}

Invariant manifolds provide geometric structures  that describe
  dynamical behavior of nonlinear     systems. Dynamical
  reductions to attracting invariant manifolds  or dynamical
  restrictions to other (not necessarily
  attracting) invariant manifolds are often sought to gain
  understanding of nonlinear dynamics.

There have been recent  works on invariant manifolds for
stochastic or random ordinary differential equations   by
Carverhill \cite{Carverhill}, Wanner \cite{Wanner2}, Arnold
\cite{Arnold}, Boxler \cite{boxler1, boxler2}, and Mohammed
\cite{Mohammed}, among others.
These authors use the (sample-wise) cocycle property for the
solution operator of the stochastic differential equations, the
Osledets' multiplicative ergodic theorem \cite{Arnold}, and a
 less-physical but technically convenient  random norm, to prove
the existence of invariant manifolds. The construction of a random
norm needs the knowledge of Oseledets spaces (a kind of eigenspace
in random linear algebra) as well as Lyapunov exponents, both are
hardly ever available; see \cite{Arnold}, p. 191 and p.379. Random
norms are not realistic in this sense, and thus representations of
invariant manifolds and the dynamical reductions are difficult to
achieve when random norms are used.

Earlier approaches on deriving dynamical reductions on stochastic
center-like manifolds by series expansions are considered by
Knoblock and Wiesenfeld \cite{Knobloch}, Schoner and Haken
\cite{Schoner}, and Xu and Roberts \cite{XuRoberts}.


 For stochastic dynamical systems, there are various concepts
 for invariance in the
 definition of invariant manifolds \cite{WaymireDuan}. In the framework
of cocycles \cite{Arnold}, the suitable concept for invariance of
a random set  is that each orbit starting inside it stays inside
it sample-wise, modulo the change of sample due to noise. Another
  concept is almost sure invariance of a deterministic set under
  stochastic dynamics, i.e., each orbit starting inside it stays
  inside it almost surely.


\bigskip

In this paper, we consider invariant manifold reductions or
restrictions for Stratonovich and Ito stochastic differential
equations in Euclidean spaces.

\medskip

We first study the system of Stratonovich stochastic differential
equations in $\R^n$:
\begin{equation}\label{system2}
dX=[AX+F^\epsilon(X)]dt+ B(X) \circ dW(t),\;\; X(0)=x_0,
\end{equation}
where $X=X(t, \om)$ is the unknown variable; $A$ is  a $n\times n$
matrix with $k$ eigenvalues of zero real parts and $n-k$
eigenvalues of negative real parts;
 $F^\epsilon:\R\rightarrow \R$ and  {$B: \R^n\rightarrow \R^{n\times n} $}
 are   nonlinear vector and matrix functions (with $\e>0$ a small parameter), respectively;
 and $W(t)$ is a standard vector Brownian
motion (or Wiener process) taking values in $\R^n$. Moreover,
$\circ$ denotes the stochastic differential in the sense of
Stratonovich.

\medskip

Then we consider the  following   stochastic system defined by
  Ito stochastic differential equations in $\R^n$:
\begin{equation} \label{system}
dX   =  F(X )d t + B(X )dW(t), \;\; X(0) =  x_0,
\end{equation}
where again $F$ and $B$ are    vector and matrix functions in
$\R^n$ and $\R^{n\times n}$, respectively. And $W(t)$ are standard
vector Brownian motion  in $\R^n$.

\medskip

Note that the Stratonovich stochastic differential $B(X) \circ
dW(t)$ and Ito stochastic differential $ B(X) dW(t)$ are
interpreted through their corresponding definitions of stochastic
integrals \cite{Oksendal}:
$$
\int_0^T B(X) \circ dW(t) := \mbox{mean-square} \lim_{\Delta t_j
\to 0} \sum_j B(X(\frac{t_{j+1}-t_j}2)) (W_{t_{j+1}}-W_{t_j}),
$$
$$
\int_0^T B(X)  dW(t) := \mbox{mean-square} \lim_{\Delta t_j \to 0}
\sum_j B(X(t_j)) (W_{t_{j+1}}-W_{t_j}).
$$
Note the difference in the sums: In Stratonovich integral, the
integrand is evaluated at the midpoint $\frac{t_{j+1}-t_j}2$  of a
subinterval $(t_j, t_{j+1})$, while for Ito integral, the
integrand is evaluated at the left end point $t_j$. See
\cite{Oksendal} for the discussion about the difference in
physical modeling by these two kinds of stochastic differential
equations. There are also dynamical differences for these two type
of stochastic equations, even at linear level \cite{CL}.

\bigskip

In this paper, we derive  an invariant manifold reduction or
restriction principle for the above systems of stochastic
differential equations.

For the Stratonovich stochastic system (\ref{system2}), we
consider random invariant center manifolds. The dynamical
reduction is achieved by investigating asymptotic behavior of
random center manifolds.

For the Ito stochastic system (\ref{system}), we study
deterministic almost sure  invariant manifolds,   which are not
necessarily   attracting. We reformulate the local invariance
condition as invariance equations, i.e., first order partial
differential equations, and then solve these equations by the
method of characteristics. Although the local invariant manifold
is deterministic, the restriction of the original stochastic
system on this deterministic local invariant manifold is still a
stochastic system but with reduced dimension.



This paper is organized as follows: In Section \ref{s2}, we recall
some basic concepts for stochastic dynamical systems. We consider
random center manifold reduction in Section \ref{s2.5}. Finally,
in Section \ref{s3}, we  construct deterministic invariant
manifolds by investigating first order partial differential
equations via the method of characteristics, and thus obtain
dynamical restrictions for  stochastic dynamic systems.


\section{Stochastic dynamical systems}\label{s2}

In this section we introduce some   definitions in stochastic
dynamical systems,  as well as recall some usual notations
  in probability.


We consider stochastic systems in the state space $\R^n$, with the
usual  metric or distance $d(x,y)=\sqrt{\sum_{j=1}^n
(x_j-y_j)^2}$, norm or length $\|x\| =\sqrt{\sum_{j=1}^n  x_j
^2}$, and the usual scalar product $<x,   y>= \sum_{j=1}^n x_j
y_j$. All invariant manifolds and their sample versions are in
this state space.

Some stochastic processes, such as a Brownian motion, can be
described by a canonical (deterministic) dynamical system (see
\cite{Arnold}, Appendix A). A standard Brownian motion (or Wiener
process) $W(t)$ in $\R^n$, with two-sided time $t \in \mathbb{R}$,
is a stochastic process with $W(0)=0$ and stationary independent
increments satisfying $W(t)-W(s) \thicksim \mathcal{N} (0,
|t-s|I)$. Here $I$ is the $n\times n$ identity matrix. The
Brownian motion can be realized in a canonical sample space of
continuous paths passing the origin at time $0$
\[
\Om = C_0(\R, \R^n): =\{\om \in C(\R,\R^n): \om(0)=0 \}.
\]
The  convergence concept in this  sample space is the  uniform
convergence on bounded and closed time intervals,  induced by the
following   metric
\[
\rho(\om, \om'):= \sum_{n=1}^{\infty}\frac1{2^n}\;
\frac{\|\om-\om'\|_n}{1+\|\om-\om'\|_n},\; \mbox{where}\;
\|\om-\om'\|_n:=\sup_{-n\leq t \leq n} \|\om(t)-\om'(t)\|.
\]
With this metric, we can define events represented by open balls
in $\Om$. For example, a ball centered at zero with radius $1$ is
$\{\om: \; \rho(\om, 0) < 1 \}$. We define the Borel
$\sigma-$algebra $\mathcal{F}$ as the collection of events
represented by    open balls $A$'s,
 complements of open balls, $A^c$'s, unions and intersections of $A$'s and/or $A^c$'s,
together with the empty event,   the whole event (the sample
 space $\Om$), and all events formed by
 doing the complements, unions and intersections forever in this collection.

 Taking the
(incomplete) Borel $\sigma-$algebra $\mathcal{F}$ on $\Om$,
together with the corresponding Wiener measure $\PX$, we obtain
the canonical probability space $(\Omega, \mathcal{F},
\mathbb{P})$, also called the Wiener space. This is similar to the
game of gambling with a dice, where the canonical sample space is
$\Om_{dice}=\{1, 2, 3, 4, 5, 6  \}$. Moreover, $ \mathbb{E}$
  denotes the mathematical expectation with respect to probability $ \mathbb{P}$.

The canonical \emph{driving} dynamical system describing the
Brownian motion is defined as
\[
\theta(t): \Om \to \Om,\;\;\; \theta(t)\om(s):=\om(t+s)-\om(t),\;
s, t \in \mathbb{R}.
\]
Then $\theta(t)$, also denoted as $\theta_t$, is a homeomorphism
for each $t$ and $(t, \om) \rightarrowtail \theta(t)\om$ is
continuous,   hence measurable. The Wiener measure $\PX$ is
invariant and ergodic under this so-called Wiener shift
$\theta_t$. In summary, $\theta_t$ satisfies the following
properties.
\begin{itemize}
    \item $\theta_0=id$,
    \item $\theta_t\theta_s=\theta_{t+s}$, $\;\;$ for all $s$,
    $t\in\R$,
    \item the map $(t,\omega)\mapsto \theta_t\omega$ is
    measurable and $\theta_t\mathbb{P}=\mathbb{P}$ for all
    $t\in\R$.
\end{itemize}

We now introduce an important concept. A filtration is an
increasing family of information accumulations, called
$\sigma$-algebras, $\mathcal{F}_t$. For each $t$, $\sigma$-algebra
$\mathcal{F}_t$ is a collection of events in sample space $\Om$.
One might observe the Wiener process $W_t$ over time $t$ and   use
  $\mathcal{F}_t$ to represent the information accumulated up to
and including time $t$.  More formally, on  $(\Omega,
\mathcal{F})$, a filtration is a family of $\sigma$-algebras
${\mathcal{F}_s : 0 \geq s \leq t}$ with $\mathcal{F}_s$ contained
in $\mathcal{F}$ for each $s$, and $ \mathcal{F}_s \subset
\mathcal{F}_{\tau}$ for $s \leq \tau$. It is also useful to think
$\mathcal{F}_t$ as the $\sigma$-algebra generated by infinite
union of $\mathcal{F}_s$'s, which is contained in $\mathcal{F}_t$.
So a filtration is often used to represent the change in the set
of events that can be measured, through gain or loss of
information.

For understanding stochastic differential equations from a
dynamical point of view, the natural filtration is defined as a
two-parameter family of  $\sigma$-algebras generated by increments
\[
\mathcal{F}_s^t:= \sigma(\om(\tau_1)-\om(\tau_2): s\leq \tau_1,
\tau_2 \leq t), \;\; s, t \in \mathbb{R}.
\]
This represents the information accumulated from time $s$ up to
and including time $t$. This   two-parameter   filtration allows
us to define forward as well as backward stochastic integrals, and
thus we can solve a stochastic differential equation from an
initial time forward as well as backward in time \cite{Arnold}.

The solution operator for the stochastic system (\ref{system2}) or
(\ref{system}) with initial condition $x(0)=x_0$ is denoted as
$\phi(t, \om, x_0)$.

The dynamics of the system on the state space $\R^n$, over the
driving flow $\theta_t$ is described by a cocycle. A cocycle
$\phi$ is a mapping:
\[
\phi:\mathbb{R} \times \Omega\times \R^n \to \R^n
\]
which  is
$(\mathcal{B}(\mathbb{R})\otimes\mathcal{F}\otimes\mathcal{B}(\R^n),\mathcal{F})$-measurable
 such that
\begin{eqnarray*}
&\phi(0,\omega,x)=x \in \R^n,\\
&
\phi(t_1+t_2,\omega,x)=\phi(t_2,\theta_{t_1}\omega,\phi(t_1,\omega,x)),
\end{eqnarray*}
for $t_1,\,t_2\in\mathbb{R},\,\omega\in \Omega,$ and $x\in \R^n$.
Then $\phi$, together with the driving dynamical system, is called
a {\em random dynamical system}. Sometimes we also use $\phi(t,
\om)$ to denote this system.

Under very general conditions,  the stochastic differential
systems (\ref{system2}) and (\ref{system}) each generates a random
dynamical system in $\R^n$; see \cite{Arnold, Imkeller}.

\medskip

We recall some concepts in   dynamical systems.  A \emph{manifold}
$M$  is a set, which locally looks like an Euclidean space.
Namely, a ``patch" of the manifold $M$ looks like a ``patch" in
$\R^n$. For example, curves, torus  and spheres in $\R^3$ are one-
and two-dimensional differentiable manifolds, respectively.
However, a manifold arising from the study of invariant sets for
dynamical systems in $\R^n$, can be very complicated. So we give a
formal definition of manifolds. For more discussions on
differentiable manifolds, see \cite{Marsden, Perko}.

\begin{definition} \textbf{(Differentiable manifold and Lipschitz manifold)}
An n-dimensional differentiable manifold $M$, is a connected
metric space with an open covering $\{U_{\alpha}\}$, i.e,
$M=\bigcup_{\alpha}U_{\alpha}$, such that

(i) for all $\alpha$ , $U_{\alpha}$ is homeomorphic to the open
unit ball in $\R^n$, $B=\{x \in \R^n :\; |x| < 1\}$, i.e., for all
$\alpha$ there exists a homeomorphism of $U_{\alpha}$ onto B,
$h_{\alpha}:U_{\alpha} \rightarrow  B$, and

(ii) if $U_{\a} \cap U_{\b}  \neq \varnothing$ and $ h_{\a}:
U_{\a} \rightarrow B$, $h_{\b}: U_{\b} \rightarrow B$ are
homeomorphisms, then $h_{\a}(U_{\a} \cap U_{\b})$ and
$h_{\b}(U_{\a}\cap U_{\b})$ are subsets of $\R^n$ and the map
\begin{equation} \label{map}
 h=h_{\a} \circ h_{\b}^{-1}: h_{\b}(U_{\a}\cap
U_{\b})  \rightarrow  h_{\a}(U_{\a}\cap U_{\b})
\end{equation}
is differentiable, and for all $x\in h_{\b}(U_{\a}\cap U_{\b})$,
the Jacobian determinant $\det Dh(x)  \neq 0$.

If the map (\ref{map}) is only Lispchitz continuous, then we call
 $M$ an n-dimensional  Lispchitz continuous manifold.
\medskip

Recall that a homeomorphism of A to B is a continuous one-to-one
map of A onto B, $h:A \rightarrow B$,  such that $h^{-1}: B
\rightarrow A$ is continues.
\end{definition}

\medskip

Just as invariant sets are important building blocks for
deterministic dynamical systems,   invariant sets are basic
geometric objects to help understand stochastic dynamics. Here we
present two different concepts about invariant sets for stochastic
systems: random invariant sets and almost sure invariant sets.

\medskip



\begin{definition} \textbf{(Random invariant set)}
A random set $M(\omega)$ is called an invariant set for a random
dynamical system $\phi$ if
\[
\phi(t,\omega,M(\omega))\subset M(\theta_t\omega), \; \;  t \in
\mathbb{R} \;\; \mbox{and}\;\; \om \in \Om.
\]
\end{definition}

\begin{definition} \textbf{(Random invariant manifold)}
If  a random invariant set  $M$ can be represented by a graph of a
Lipschitz mapping
\[
\gamma^\ast(\omega,\cdot): \;  H^+\to H^-, \; \mbox{with direct
sum decomposition} \quad H^+\oplus H^-= \R^n
\]
such that
\begin{eqnarray*}
M(\omega)=\{x^++\gamma^\ast(\omega,x^+),x^+\in H^+\},
\end{eqnarray*}
then $M$ is called a Lipschitz continuous invariant manifold.
\end{definition}




\bigskip

We will also consider   deterministic invariant sets or manifolds,
while the invariance   is  in the sense of almost-sure (a.s.)
\cite{Aubin, Filipovic, Chueshov, Stanzhitskii, Zab}.

\begin{definition} \textbf{(Almost sure invariant set)}
A (deterministic) set $M$ in $\mathbb{R}^n$ is called locally
almost surely invariant for (\ref{system}), if for all $(t_0,x_0)
\in \R_+ \times M $, there exists a continuous local weak solution
$X^{(t_0,x_0)}$ with lifetime $\tau=\tau(t_0,x_0)$, such that
\begin{eqnarray*}
X^{(t_0,x_0)}_{t\wedge \tau} \in M , \;\; \forall  t > t_0, \;\;\;
a.s. \;\;  \om \in \Om,
\end{eqnarray*}
where $t\wedge \tau = \min(t, \tau)$.
\end{definition}


\section{Random center manifold reduction}  \label{s2.5}

In this section we study the system of Stratonovich stochastic
differential equations in $\R^n$:
\begin{equation}  \label{sode}
dX=[AX+F^\epsilon(X)]dt+B(X) \circ dW(t),\;\; X(0)=x_0,
\end{equation}
where  $A$ is  a $n\times n$ matrix with $k$ eigenvalues of zero
real parts and $n-k$ eigenvalues of negative real parts;
 $F^\epsilon:\R\rightarrow \R$ and  {$B: \R^n\rightarrow \R^{n\times n} $}
 are   Lipschitz continuous  with Lipschitz constants $L_F^\epsilon$ and
 $L_B$, respectively;
$F^\epsilon(0)=B(0)=0$.
  And
$\epsilon$ is a small parameter so that $F^\epsilon$ can be seen
as a small perturbation, that is, we have $L_F^\epsilon\rightarrow
0$ as $\epsilon\rightarrow 0$. Let $X_c$ be the projection of $X
\in \R^n$ on $\ker A$ (i.e., projection to the center eigenspace).
Note that later on we will truncate the nonlinearity so that it
has global Lipschitz constant. For stochastic systems, truncation
may not be always appropriate, although sometimes it works fine,
such as in considering nonlinear dynamical behavior near fixed
points \cite{CLR}.  We have the following result.

\medskip

\begin{theorem} \label{center}(\textbf{Random center manifold reduction})

 Given the above  assumptions for   the  Stratonovich stochastic
differential equation   (\ref{sode}). If further assume that the
equation  (\ref{sode})  generates a dissipative random dynamical
system (e.g., has an absorbing set). Then for small $\epsilon$,
the long time behavior of (\ref{sode}) can be described by the
following stochastic differential equation in $\R^k$
\begin{equation}\label{e2}
dX_c(t)=F^\epsilon(X_c(t))dt+ B(X_c(t)) \circ dW_c(t)
\end{equation}
provided (\ref{e2}) is structurally stable. In this reduced
system, $F^\epsilon(X_c):=F^\epsilon(X_c+0)$ and
$B(X_c):=B(X_c+0)$ with $0\in \R^n \backslash \ker A$,
and $W_c(t)$ is the projection of $W(t)$ on $\ker A$.  \\
\end{theorem}

 \textbf{Remark 1}   We say the long time dynamics of the stochastic equation
 (\ref{sode}) is described by the stochastic equation (\ref{e2}) if
  both systems have the same limit sets
  (and possibly also share some other invariant sets).

\medskip

 \textbf{Remark 2}  The random dynamical system $\phi(t,\om)$
generated by (\ref{e2}) is called structurally stable, if for any
small perturbation (small in the sense of the usual metric in the
space of continuous functions) to $F^\epsilon(x)$ and $B(x)$, the
perturbed random dynamical system $\Phi(t,\om)$ is topologically
equivalent to $\phi(t, \om)$. Namely, there exists a random
homeomorphism $h(\om)$ so that $\Phi(t,\om) \circ h(\om) =
h(\theta_t \om) \circ \phi(t, \om)$.

\medskip

 The proof for this theorem can be obtained by modifying  the proof in \cite{WangDuan}
 to the systems in the space   $\R^n$.

Let us look at an   example.

\begin{example}
Consider a system of stochastic  differential equations in
$\mathbb{R}^2$:
\begin{eqnarray*}
dx &=&  -xdt +   (xy^2 -x^3-\frac12 x)dt + x \circ dW_1(t), \\
dy &=&  0 ydt +   (-2+x^2y -y^3 -\frac12 y)dt + y \circ dW_2(t),
\end{eqnarray*}
where $W_1$ and $W_2$ are independent scalar Brownian motions.
\end{example}

Let $u=(x, y)^T$. Then
\begin{eqnarray*}
du =  (Au + \tilde{F}(u))dt +  Bu \circ dW(t),
\end{eqnarray*}
with
\begin{eqnarray*}
A = \left( \begin{array}{ll}
         -1 & 0\\
          0 & 0\end{array} \right),
\tilde{F}(u) = \left( \begin{array}{ll}
         xy^2 - x^3-\frac12 x\\
         -2+x^2y - y^3-\frac12 y\end{array} \right)  \; \mbox{and} \;
 Bu \circ dW(t) = \left( \begin{array}{ll}
         x \circ dW_1(t)\\
         y \circ dW_2(t)\end{array} \right).
\end{eqnarray*}

In order to apply the Ito's formula, we rewrite this system in the
equivalent Ito's stochastic differential equations (see
\cite{Oksendal}, page 36):
\begin{eqnarray*}
dx &=&  -xdt +   (xy^2 -x^3)dt + x dW_1(t), \\
dy &=&  0 ydt +  (-2+x^2y -y^3)dt + y dW_2(t),
\end{eqnarray*}
where $W_1$ and $W_2$ are independent scalar Brownian motions. Let
$u=(x, y)^T$, then
$$
du=(Au+F(u))dt+Bu\; dW(t)
$$
with
$$ A=\left(
  \begin{array}{c l}
    -1 & 0\\
    0 & 0
    \end{array}
\right), F(u)=\left(
  \begin{array}{c l}
    xy^2-x^3 \\
    -2+x^2y-y^3
    \end{array}
\right)\;\; \mbox{and} \;\;  Bu=\left(
  \begin{array}{c l}
    x & 0\\
    0 & y
    \end{array}
\right).
 $$
 Recall the standard scalar product $<u_1, u_2>=x_1x_2+y_1y_2$ and norm
 $\|u\|=\sqrt{x^2+y^2}$ in $\mathbb{R}^2$.
Then, we apply the Ito's formula   (see \cite{Oksendal}, page 48)
to obtain ``energy" estimate
\begin{eqnarray*}
\frac{1}{2}\frac{d}{dt}\mathbb{E}\|u\|^2&=&\mathbb{E}\langle u,
du\rangle + \mathbb{E}{1 \over 2}
\langle du, du\rangle \\
&=&\mathbb{E}\langle u, du\rangle  + \mathbb{E}{1 \over 2} \langle Bu\; dW(t), Bu\; dW(t)\rangle \\
&=&\mathbb{E}\langle u, A(u)+F(u)\rangle  + {1 \over 2}\mathbb{E} \; Trace[Bu \cdot(Bu)^T] \\
&=& -x^2 + x^2y^2 - x^4 -2y + x^2y^2 -y^4+{1 \over 2} (x^2 + y^2)\\
&=& -{1 \over 2}(x^2+y^2) + (2x^2y^2 - x^4 -y^4) - y + y^2 -1 + 1 \\
&=& -{1 \over 2}(x^2+y^2) - (x^2-y^2)^2 + (y-1)^2 - 1\\
 &\leq& -{1 \over 2} E\|u\|^2,
\end{eqnarray*}
if $ y$ is near the equilibrium point $(0, 0)$ (so that $0<y<1$).
Note that here $ \mathbb{E}$
  denotes the expectation with respect to probability $ \mathbb{P}$.
This estimate will be used to conclude dissipativity for the
(truncated) system.

The nonlinear terms $(xy^2 -x^3)$  and $(-2+ x^2y -y^3)$ can be
truncated within a disk centered at $(0, 0)$ with radius $0<\e <
1$ (making them zero outside the disk). The truncated nonlinear
terms satisfy desired Lipschitz conditions.
  And, the above   ``energy"  estimate   implies the dissipative property for the truncated system.


By  Theorem \ref{center}, near the equilibrium  point $(0,0)$
(i.e., taking $\e$ is small enough),  the  original
two-dimensional  system is asymptotically reduced to a
one-dimensional stochastic dynamical system
\begin{eqnarray*}
dy  =    (-2 -   y^3) dt +y dW_2(t).
\end{eqnarray*}


\section{Invariant manifold restriction }\label{s3}

Now we consider the      stochastic system (\ref{system}) defined
by
  Ito stochastic differential equations in $\R^n$:
\begin{equation} \label{Ito}
dX   =  F(X )d t + B(X )dW(t), \;\; X(0) =  x_0,
\end{equation}
where   $F$ and $B$ are    vector and matrix functions in $\R^n$
and $\R^{n\times n}$, respectively.   We also assume that
  $F( \cdot) \in
C^1(\R^n;\R^n) $ and $B( \cdot) \in C^1(\R^n;\R^{n \times n})$.

\bigskip



  We are going to derive representations of   invariant finite dimensional
manifolds in terms of $A , F$ and $B $,  by using the tangency
conditions for a deterministic $C^2$ manifold $M$ in $\R^n$:
\begin{eqnarray}
\mu(\om,x) := F(\om, x) &-& {1 \over
2}\sum_j[DB^j(\om, x)]B^j(\om, x) \in T_xM,    \label{c1}    \\
B^j(\om, x) &\in& T_xM , \;\;    j=1, \cdots, n,  \label{c2}
\end{eqnarray}
where $D$ represents Jacobian operator and $B_j$ is the $j-$th
column of the matrix $B$. The above   tangency conditions are
shown to be equivalent to   almost sure   local invariance of
manifold $M$; see Filipovic (\cite{Filipovic}) and   related works
\cite{Aubin, Milian, Zab, Chueshov, Amann}.

The almost sure invariance conditions (\ref{c1})-(\ref{c2}) for
manifold $M$ mean that the $n+1$ vectors, $\mu$ and $B^j, \;  j=1,
\cdots, n$, are tangent vectors to $M$. Namely, these $n+1$
vectors are orthogonal to the normal vectors of manifold $M$.

In other words, if the normal vector  for $M$ at $x$ is $N(x)$,
then the almost sure invariance conditions (\ref{c1})-(\ref{c2})
become the following \emph{invariance equations} for manifold $M$:
For all $x\in M$,
\begin{eqnarray}
 \mu (x) \cdot N(x) & =& 0,    \label{c3}     \\
 B^j(x) \cdot N(x) &= & 0, \;\;    j=1, \cdots, n,  \label{c4}
\end{eqnarray}
where $\cdot$ denotes the scalar product in $\mathbb{R}^n$.

Invariant manifolds are usually represented as graphs of some
functions in $\R^n$. By investigating the above invariance
equations (\ref{c3})-(\ref{c4}), we may be able to find some local
invariant manifolds $M$ for the stochastic system (\ref{Ito}).

The goal for this section is to present a method to find some of
these local invariant manifolds. Although  the following result
and example are stated for a codimension $1$ local invariant
manifold, the idea extends to other lower dimensional local
invariant manifolds, as long as the normal vectors $N(x)$ (or
tangent vectors) may be represented; see tangency conditions
(\ref{c3})-(\ref{c4}) above and (\ref{c5})-(\ref{c6}) below.


\begin{theorem} \label{local}
(\textbf{Local invariant manifold restriction}) \\
 Let the local invariant manifold $M$ for the stochastic dynamical system (\ref{Ito}) be represented as a
graph defined by the algebraic equation
\begin{eqnarray} \label{G}
M: \;\;\; G(x_1 , \cdots , x_n) = 0.
\end{eqnarray}
Then $G$ satisfies a system of first order (deterministic) partial
 differential equations and the local invariant manifold $M$ may
 be found by solving these partial
 differential equations by the method of characteristics. By
 restricting the original dynamical system (\ref{Ito}) on this
 local invariant manifold $M$, we obtain a locally valid, reduced lower
 dimensional system.
 \end{theorem}

In fact, the normal vector to this graph or surface is, in terms
of partial derivatives, $\nabla G(x) =( G_{x_1} , \cdots , G_{x_n}
) $. Thus the invariance equations (\ref{c3})-(\ref{c4}) are now
\begin{eqnarray}
 \mu (x) \cdot  \nabla G(x) & =& 0,    \label{c5}     \\
 B^j(x) \cdot \nabla G(x) &= & 0, \;\;    j=1, \cdots, n,  \label{c6}
\end{eqnarray}
This is a system of first order partial differential equations in
$G$. We apply the method of characteristics to solve for $G$, and
therefore obtain the invariant manifold $M$, represented by a
graph in state space $\mathbb{R}^n$:
 $G(x_1 , \cdots , x_n) = 0$.\\

\medskip

In the rest of this section, we first recall the method of
characteristics, and then work out an example of finding a local
invariant manifold and reduced system.

\medskip

 {\bf Method of
Characteristics}: Consider a first order partial differential
equation for the unknown scalar function $u$ of n variables $x_1 ,
\cdots , x_n$
\begin{eqnarray} \label{first}
\sum_{j=1}^n a_i(x_1, \cdots, x_n) u_{x_i} = c(x_1, \cdots, x_n),
\end{eqnarray}
with continuous coefficients $ a_i$'s and $c$.

Note that the solution  surface $u=u(x_1,..,x_n,t)$ in
$x_1...\cdots x_nu-$space has normal vectors $N: =(u_{x_1},
\cdots, u_{x_n}, -1)$. This partial differential equation implies
that the   vector $ V= :(a_1, \cdots, a_n, c)$  is perpendicular
to this normal vector and hence must lie in the tangent plane to
the graph of $z=u(x_1, \cdots, x_n)$.

In other words, $ (a_1, \cdots, a_n, c)$  defines a vector field
in $\mathbb{R}^n$, to which graphs of the solutions must be
tangent at each point \cite{McOwen}. Surfaces that are tangent at
each point to a vector field in $\mathbb{R}^n$ are called {\bf
integral surfaces } of the vector field. Thus to find a solution
of equation (\ref{first}), we should try to find integral
surfaces.

How can we construct integral surfaces? We can try using the
characteristics curves that are the integral curves of the vector
field. That is, $X=(x_1(t), \cdots, x_n(t) )$ is a {\bf
characteristic} if it satisfies the following system of ordinary
differential equations:
\begin{eqnarray*}
{dx_1 \over dt} & = & a_1(x_1,\cdots,x_n), \\
& \cdots &\\
{dx_n \over dt}& = & a_n(x_1,..,x_n),   \\
{du \over dt} & = & c(x_1,..,x_n).
\end{eqnarray*}
A smooth union of characteristic curves is an integral surface.
there may be many integral surfaces. Usually an integral surface
is determined by requiring it to contain (or pass through) a given
initial curve  or an $n-1$ dimensional manifold $\Gamma$:
\begin{eqnarray*}
x_i  & = & f_i(s_1,..,s_{n-1}),i=1..n\\
u  & = & h(s_1,..,s_{n-1})
\end{eqnarray*}
This generates an $n$-dimensional integral manifold $M$
parameterized by $(s_1,..,s_{n-1}, t)$. The solution $u(x_1,
\cdots, x_n)$ is obtained by solving for $(s_1,..,s_{n-1},t)$ in
terms of variables $(x_1, \cdots, x_n)$.


{\bf Remark}: If initial data $\Gamma$ is non-characteristic,
i.e., it is nowhere tangent to the vector field $V=(a_1, \cdots,
a_n, c)$, and $a_1, \cdots, a_n, c$ are $C^1$ (and thus locally
Lipschitz continuous), then there exists   a unique integral
surface $u=u(x_1, \cdots, x_n)$ containing $\Gamma$, defined at
least locally near $\Gamma$.

\bigskip

Now applying the above method of characteristics to
(\ref{c5})-(\ref{c6}), we obtain a solution $G=G(x_1, \cdots,
x_n)$.  However, the local invariant manifold  $M$ that we look
for is represented by the equation
$$
G(x_1, \cdots, x_n) =0.
$$
Therefore, a skill is needed to make sure that the solution
$G=G(x_1, \cdots, x_n)$ actually penetrates the plane $G=0$ in the
$x_1\cdots x_n G-$space; see Fig. 1. This needs to be achieved by
slecting appropriate initial data $\Gamma$.  The invariant
manifold $M$ we thus obtain is defined at least locally near the
initial data $\Gamma$.


\bigskip

We illustrate the method for finding local invariant manifold and
the corresponding   reduced system by an example.

\begin{example}

\begin{eqnarray*}
{dx \over dt} & = &  x  + x \;\dot W_1  + x \; \dot W_2,\\
{dy \over dt} & = & 3x+2y + (x+y) \;\dot W_1 + (x+y) \;\dot W_2
\end{eqnarray*}
where $W^1_t$ and $W^2_t$ are independent scalar Brownian motions.
\end{example}

We look for a local invariant manifold $M \subset \mathbb{R}^2$.
For this illustrative example, the associated tangency conditions
(\ref{c1}) and (\ref{c2}) coincide and thus becomes a single
\emph{invariance condition}:
\begin{eqnarray} \label{inv}
(x,x+y)^T &\in& T_xM
\end{eqnarray}

We represent the invariant manifold $M$ by $G(x,y) = 0$. This
surface has normal vector $(G_x,G_y)$. By noticing that normal
vector is orthogonal to the tangent surface $T_xM$, we see that
the above single invariance condition (\ref{inv}) becomes a single
invariance equation:
\begin{eqnarray*}
xG_x + (x+y)G_y &=& 0.
\end{eqnarray*}
We    solve this first order partial differential equation with
initial curve $\Gamma$ parameterized as $(f(s), g(s), h(s))$. The
characteristic equations are
\begin{eqnarray*}
{dx \over dt} & = & x ,\\
{dy \over dt} & = & x+y ,\\
{dG \over dt} & = & 0.
\end{eqnarray*}
We   solve these equations and invoke the initial conditions to
find that
\begin{eqnarray*}
x & = & f(s)e^t ,\\
y & = & (f(s)t+g(s))e^t ,\\
G & = & h(s).
\end{eqnarray*}
This is the general solution with respect to the general initial
condition $(x_0(s), y_0(s), G_0(s)): =(f (s), g(s),h(s))$.  By
solving for $t, s$ in terms of $x, y$, we obtain $G=G(x,y)$.

We illustrate this by a specific choice of initial curve $(f(s),
g(s), h(s))$.  Note that, in order to obtain a local invariant
manifold $G(x, y)=0$, we also need to pick initial curve so that
$G$ actually takes both positive and negative values, and thus the
invariant manifold $G(x,y)=0$ is defined on some set in the
$xy-$plane. In other words, the continues function $G(x,y)$
  satisfies $max\{G\}*min\{G\} \leq 0$ locally.

For example, taking $\Gamma: (x_0(s), y_0(s), G_0(s)) =(1, s, s)$,
we  then have
\begin{eqnarray*}
x & = & e^t,\\
y & = & (t+s) e^t ,\\
G & = & s.
\end{eqnarray*}
Thus $s = {y \over x}-ln(x)$ and $G(x,y)={y \over x}-ln(x)$. Thus
an invariant manifold $M$ is $G(x,y)=0$, i.e.,
\begin{eqnarray*}
   {y \over x}-ln(x)  =0.
\end{eqnarray*}



\bigskip
\bigskip

{\bf Acknowledgements.}

This work was partly supported by the NSF Grants DMS-0209326  \&
DMS-0542450. A part of this work was done while J. Duan was
visiting the American Institute of Mathematics, Palo Alto,
California, USA.
  We thank Hans Crauel, Xinchu Fu, Kening Lu,   Bjorn
Schmalfuss and Wei Wang for helpful discussions.



\begin{thebibliography}{50}

\bibitem{Marsden} R. Abraham, J. E. Marsden and T. Ratiu.
\emph{Manifolds, Tensor Analysis, and Applications}. Second Ed.,
Springer-verlag, New York, 1988.

\bibitem{Amann} H. Amann.
Invariant sets and existence theorems for semilinear parabolic and
elliptic systems. \emph{J. Math. Anal. Appl.}  \textbf{65} (1978),
432-467.

\bibitem{Arnold}
L. Arnold.
\newblock {\em Random Dynamical Systems.}
\newblock Springer-Verlag, New York, 1998.


\bibitem{Aubin}
J.-P. Aubin and G. Da Prato.
\newblock Stochastic Viability and Invariance.
\newblock {\em  Scuola Norm. Sup. Pisa}, {\bf l27} (1990), 595-694.


\bibitem{boxler1}
P. Boxler.
\newblock How to constract stochastic center manifolds on the level of vector fields.
\newblock {\em  Lectures Notes in Math}, {\bf 1486} (1991), 141-158.


\bibitem{boxler2}
P. Boxler.
\newblock A stochastic version of center manifold theory.
\newblock {\em  Probability Theory \& Related Fields}, {\bf 83} (1989), 509-545.






\bibitem{CLR}
T. Caraballo,  J. Langa   and J. C. Robinson.
  A stochastic pitchfork bifurcation in a reaction-diffusion equation.
 {\em R. Soc. Lond. Proc. Ser. AMath. Phys. Eng. Sci.},
{\bf 457}(2001), no.2013 2041-2061.


\bibitem{CL}
T. Caraballo and J. Langa.
\newblock A comparison of the longtime
behavior of linear Ito and Stratonovich partial differential
equations.
\newblock {\em Stochastic Anal. Appl.}, {\bf 19}(2001), no.2 183-195.



\bibitem{Carverhill}
A. Carverhill.
\newblock Flows of stochatsic dynamical systems: Ergodic theory.
\newblock \emph{Stochastics}, {\bf 14} (1985), 273-317.

\bibitem{Chueshov}
I. D. Chueshov and P.-A. Vuillermot.
\newblock Non-random Invariant Sets for Some systems of Parabolic
Stochatstic Partial Differentia Equations.
\newblock \emph{Stochastic Analysis and Apllication}, Vol.22 , No.6,pp1421-1486, 2004.

\bibitem{DuanLuSchm}
J. Duan, K. Lu  and B. Schmalfu{\ss}. Invariant manifolds for
stochastic partial differential equations.
  {\em The Annals of Probability}, {\bf 31}(2003), 2109-2135.


\bibitem{DuanLuSchm2} J. Duan, K. Lu and B.  Schmalfu{\ss}.
Smooth stable and unstable manifolds for stochastic evolutionary
equations. {\em J. Dynamics and Diff. Eqns.} {\bf 16} (2004),
949-972.


\bibitem{Filipovic}
\newblock D. Filipovic.
\newblock Invariant manifolds for weak solutions to  stochastic equations.
\newblock {\em Probability Theory \& Related Fields }, Volume \textbf{118} (2000),
Number 3. 323 - 341.


\bibitem{Gas} V. A. Gasanenko.
On invariant sets for stochastic differential equations.
\emph{Theory of Stochastic Processes} \textbf{9} (2003), 60-64.


\bibitem{GH}
J. Guckenheimer and P. Holmes.
\newblock \emph{Nonlinear Oscillations,Dynamical Systems and Bifurcations
of Vector Fields}.
\newblock Springer-Verlag, New York, 1983.


\bibitem{Holden} H. Holden, B. Oksendal, J. Uboe and T. Zhang.
\emph{Stochastic Partial Differential Equations}. Birkhauser,
Boston, 1996.


\bibitem{Imkeller} P. Imkeller and C. Lederer.   The cohomology of stochastic and random differential equations,
 and local linearization of stochastic flows. \emph{Stoch. Dyn.} \textbf{2}
 (2002), no. 2, 131--159.


\bibitem{Knobloch} E. Knobloch and K. A. Wiesenfeld.
Bifurcations in fluctuating systems: The center-manifold approach.
\emph{J. Stat. Phys.} \textbf{33} (1983), 611-637.


\bibitem{Kunita}
H. Kunita.
\newblock {\em Stochastic flows and stochastic differential
equations.}
\newblock  Cambridge University Press, 1990.


\bibitem{Milian}
A. Milian.
\newblock Invariance for stochastic equations with regular
coeffiencients,
\newblock {\em Stochastic Anal. Appl.}, {\bf 15}(1997), no.1. 91-101.


\bibitem{Mohammed}
S.-E.A. Mohammed and M. Scheutzow.
\newblock The Stable Manifold Theorem for Stochastic Differential
Equations.
\newblock {\em The Annals of  Probability}, Vol. \textbf{27}, No. 2,
(1999), 615-652.

\bibitem{McOwen}
R.C McOwen.
\newblock \emph{Partial Differential Equantions}.
\newblock Pearson Education, New Jetsy, 2003.

\bibitem{Oksendal}
B. Oksendal.
\newblock {\em Stochastic Differenntial Equations.} Sixth Ed.,
\newblock Springer-Verlag, New York, 2003.



\bibitem{Perko}
 L. Perko.
 \newblock {\em Differential Equations and Dynamical
Systems.}
\newblock Cambridge University Press, 1990.




\bibitem{Rue79}
D.~Ruelle,
\newblock Ergodic theory of differentiable dynamical systems.
\newblock {\em Publ. Math. Inst. Hautes Etud. Sci.}   (1979), 275-306.

\bibitem{Rue82}
D.~Ruelle,
\newblock Characteristic exponents and invariant manifolds in Hilbert spaces.
\newblock {\em Ann. of Math.}  {\bf 115} (1982), 243--290.



\bibitem{Schoner} G. Schoner and H. Haken.
A systematic elimination procedure for Ito stochatsic differential
equations and the adiabatic approximation. \emph{Z. Phys. B}
\textbf{68} (1987), 89-103.




\bibitem{Stanzhitskii} A. M. Stanzhitskii.
Investigation of invariant sets of Ito's stochastic systems with
the use of Lyapunov functions. \emph{Ukrainian Math. J. }
\textbf{53} (2001), No. 2, 323-327.


\bibitem{WangDuan}  W. Wang and J. Duan.
Invariant manifold reduction and bifurcation
       for stochastic partial differential equations.
       \emph{Submitted}, 2005.

\bibitem{Wanner2} T. Wanner.
\newblock  Linearization of Random Dynamical Systems.
\newblock {\em Dynamics Report}, Volumn 4. Spring-Verlog, New York, 1995.


\bibitem{WaymireDuan} E. Waymire  and J. Duan (Eds.).
\emph{Probability and Partial Differential Equations in Modern
Applied Mathematics}. Springer-Verlag,   2005.


\bibitem{Wiggins} S. Wiggins.
\newblock \emph{Introduction to Applied Nonlinear Dynamical Systems and Chaos}.
\newblock Spring-Verlag, New York, 1990.

\bibitem{XuRoberts}
C. Xu and A. J. Roberts,
\newblock On the low-dimensional modelling of Stratonovich stochastic
differential equations.
\newblock {\em Physica A}, {\bf 225}(1996), 62-68.

\bibitem{Zab}   J. Zabczyk. Stochastic invariance and consistency of financial
models. \emph{Atti Accad. Naz. Lincei Cl. Sci. Fis. Mat. Natur.
Rend. Lincei} (9) Mat. Appl. \textbf{11} (2000), no. 2, 67--80.








\end{thebibliography}
\end{document}